\documentclass[a4paper,11pt]{amsart}
\usepackage{graphicx}
\usepackage{mathptmx}
\usepackage{amsmath}
\usepackage{amssymb}
\usepackage{enumitem}
\usepackage{xcolor}
\usepackage{xparse}
\NewDocumentCommand{\eulerian}{omm}
 {%
  \genfrac<>{0pt}{}{#2}{#3}%
  \IfValueT{#1}{_{\!#1}}%
 }

 \newmuskip\pFqmuskip

\newcommand*\pFq[6][8]{%
  \begingroup 
  \pFqmuskip=#1mu\relax
  \mathchardef\normalcomma=\mathcode`,
  \mathcode`\,=\string"8000
  \begingroup\lccode`\~=`\,
  \lowercase{\endgroup\let~}\pFqcomma
  {}_{#2}F_{#3}{\left(\genfrac..{0pt}{}{#4}{#5}\bigg|#6\right)}%
  \endgroup
}
\newcommand{\pFqcomma}{{\normalcomma}\mskip\pFqmuskip}

\newtheorem{theorem}{Theorem}

\newtheorem{corollary}[theorem]{Corollary}

\begin{document}

\title[Some identities involving degenerate Stirling numbers]{Some identities involving degenerate Stirling numbers associated with several degenerate polynomials and numbers}

\author{Taekyun  Kim}
\address{Department of Mathematics, Kwangwoon University, Seoul 139-701, Republic of Korea}
\email{tkkim@kw.ac.kr}

\author{DAE SAN KIM}
\address{Department of Mathematics, Sogang University, Seoul 121-742, Republic of Korea}
\email{dskim@sogang.ac.kr}


\maketitle

\begin{abstract}
The aim of this paper is to investigate some properties, recurrence relations and identities involving degenerate Stirling numbers of both kinds associated with degenerate hyperharmonic numbers and also with degenerate Bernolli, degenerate Euler, degenerate Bell and degenerate Fubini polynomials. 
\end{abstract}

\section{Introduction}
It is noteworthy that various degenerate versions of quite a few special polynomials and numbers have been investigated recently, which began from the pioneering work of Carlitz in [2,3]. These explorations for degenerate versons are not only limited to special polynomials and numbers but also extended to some transcendental functions, like gamma functions. In the course of this quest, many different tools are used, which include generating functions, combinatorial methods, $p$-adic analysis, umbral calculus, operator theory, differential equations, special functions, probability theory and analytic number theory (see [8-14,16] and the references therein).\par
The aim of this paper is to investigate some properties, recurrence relations and identities involving degenerate Stirling numbers (see (10), (11)) associated with degenerate hyperharmonic numbers (see (5), (19)) and also with degenerate Bernoulli, degenerate Euler, degenerate Bell and degenerate Fubini polynomials (see (15), (16), (13), (17)). The novelty of the present paper is that a degenerate version of the hyperharmonic numbers, namely the degenerate hyperharmonic numbers, is firstly introduced. \par
The outline of this paper is as follows. In Section 1, we recall degenerate exponentials, degenerate logarithms, degenerate Stirling numbers of both kinds and hyperharmonic numbers. Also, we remind the reader of degenerate Bell polynomials, degenerate Bernoulli polynomials, degenerate Euler polynomials, degenerate Fubini polynomials and degenerate polylogarithms, which are, respectively, degenerate versions of Bell polynomials, Bernoulli polynomials, Euler polynomials, Fubini polynomials and polylogarithms. Section 2 is the main result of this paper. We introduce degenerate hyperharmonic numbers, which are a degenerate version of the hyperharmonic numbers, and derive the generating function of them in Theorem 1. We obtain an identity involving the degenerate hyperharmonic numbers and the degenerate Stirling numbers of the second kind in Theorem 2, and its inversion identity in Theorem 3. Obtained are identities involving the degenerate Euler polynomials and the degenerate Stirling numbers of both kinds in Theorem 4. Derived are identities relating the degenerate Bernoulli polynomials and the degenerate Stirling numbers of both kinds in Theorems 5, 6, 7 and 10. In Theorems 8 and 9, we get identities connecting the degenerate hyperharmonic numbers, the degenerate Bernoulli polynomials and the degenerate Stirling numbers of both kinds. A recurrence relation for the degenerate Bell polynomials are deduced in Theorem 11. An identity connecting the degenerate Stirling numbers of the second kind, the degenerate Bernoulli polynomials and the degenerate Bell polynomials is derived in Theorem 12. An explicit expression and a recurrence relation for the degenerate Fubini polynomials are obtained respectively in Theorem 13 and Theorem 14. Finally, an interesting identity on a finite sum of the degenerate Stirling numbers of the second kind is derived in Theorem 15.  \par
For any $\lambda\in\mathbb{R}$, the degenerate exponential is defined by 
\begin{equation}
e_{\lambda}^{x}(t)=\sum_{k=0}^{\infty}\frac{(x)_{k,\lambda}}{k!}t^{k},\quad (\mathrm{see}\ [9,10,13]),\label{1}	
\end{equation}
where
\begin{equation}
(x)_{0,\lambda}=1,\quad (x)_{k,\lambda}=x(x-\lambda)(x-2\lambda)\cdots(x-(k-1)\lambda),\quad (k\ge 1).\label{2}
\end{equation}
When $x=1$, we write $e_{\lambda}=e_{\lambda}^{1}(t)$. \par 
Let $\log_{\lambda}(t)$, called the degenerate logarithm, be the compositional inverse of $e_{\lambda}(t)$ such that $\log_{\lambda}(e_{\lambda}(t))=e_{\lambda}(\log_{\lambda}(t))=t$. Then we have 
\begin{equation}
\log_{\lambda}(1+t)=\sum_{k=1}^{\infty}\lambda^{k-1}(1)_{k,\frac{1}{\lambda}}\frac{t^{k}}{k!},\quad (\mathrm{see}\ [8,14]). \label{3}	
\end{equation}
Note that $\displaystyle\lim_{\lambda\rightarrow 0}e_{\lambda}^{x}(t)=e^{xt},\quad \lim_{\lambda\rightarrow 0}\log_{\lambda}t=\log t$. \par 
It is well known that the harmonic numbers are defined by 
\begin{equation}
H_{0}=0,\quad H_{n}=1+\frac{1}{2}+\frac{1}{3}+\cdots+\frac{1}{n},\quad (n\in\mathbb{N}),\quad (\mathrm{see}\ [4,5,15]).\label{4}	
\end{equation}
The generating function of the harmonic numbers is given by 
\begin{equation*}
-\frac{1}{1-t}\log(1-t)=\sum_{n=0}^{\infty}H_{n}t^{n},\quad (\mathrm{see}\ [4,5,15]).	
\end{equation*}
Recently, the degenerate harmonic numbers are introduced by Kim-Kim as
\begin{equation}
H_{0,\lambda}=0,\quad H_{n,\lambda}=\sum_{k=1}^{n}\frac{1}{\lambda}\binom{\lambda}{k}(-1)^{k-1},\quad (n\in\mathbb{N}),\quad (\mathrm{see}\ [9]). \label{5}
\end{equation}
Note that $\displaystyle\lim_{\lambda\rightarrow 0}H_{n,\lambda}=H_{n}\displaystyle$. \par 
In [5], Conway and Guy introduced the hyperharmonic numbers given by 
\begin{equation}
H_{n}^{(1)}=H_{n},\quad H_{n}^{(r)}=\sum_{k=1}^{n}H_{k}^{(r-1)},\quad (r\ge 2).\label{6}
\end{equation}
From \eqref{6}, we see that 
\begin{equation}
H_{n}^{(r)}=\binom{n+r-1}{r-1}(H_{n+r-1}-H_{r-1}),\quad (r\ge 2),\quad H_{0}^{(r)}=0,\quad (\mathrm{see}\ [5]).\label{7}
\end{equation}
The Fubini polynoimals are defined by 
\begin{equation}
\frac{1}{1-x(e^{t}-1)}=\sum_{n=0}^{\infty}F_{n}(x)\frac{t^{n}}{n!},\quad (\mathrm{see}\ [4,6,12,16]).\label{8}	
\end{equation}
By \eqref{8}, we easily get 
\begin{equation}
F_{n}(x)=\sum_{k=0}^{n}S_{2}(n,k)k!x^{k},\quad (n\ge 0),\quad (\mathrm{see}\ [12]),\label{9}
\end{equation}
where $S_{2}(n,k)$ are the Stirling numbers of the second kind given by 
\begin{displaymath}
\frac{1}{k!}(e^{t}-1)^{k}=\sum_{n=k}^{\infty}S_{2}(n,k)\frac{t^{n}}{n!},\quad (k\ge 0),\quad (\mathrm{see}\ [1,11,17,18]). 
\end{displaymath}
Recently, Kim-Kim introduced the degenerate Stirling numbers of the first kind defined by
\begin{equation}
(x)_{n}=\sum_{k=0}^{n}S_{1,\lambda}(n,k)(x)_{k,\lambda},\quad (n\ge 0),\quad (\mathrm{see}\ [8]),\label{10}
\end{equation}
where 
\begin{displaymath}
(x)_{0}=1,\quad (x)_{n}=x(x-1)(x-2)\cdots(x-n+1),\quad (n\ge 1),\quad (\mathrm{see}\ [15]). 
\end{displaymath}
As the inversion formula of \eqref{10}, the degenerate Stirling numbers of the second kind are given by 
\begin{equation}
(x)_{n,\lambda}=\sum_{k=0}^{n}S_{2,\lambda}(n,k)(x)_{k},\quad (n\ge 0),\quad (\mathrm{see}\ [8]).\label{11}
\end{equation}
From \eqref{11}, we note hat 
\begin{equation}
S_{2,\lambda}(n,k)=S_{2,\lambda}(n-1,k-1)+(k-(n-1)\lambda)S_{2,\lambda}(n-1,k),\quad (\mathrm{see}\ [8]),\label{12}	
\end{equation}
where $n,k\in\mathbb{N}$ and with $n\ge k$. \par 
In [10], the degenerate Bell polynomials are defined by
\begin{equation}
e^{x(e_{\lambda}(t)-1)}=\sum_{n=0}^{\infty}\phi_{n,\lambda}(x)\frac{t^{n}}{n!},\quad (\mathrm{see}\ [10,14]).\label{13}	
\end{equation}
When $x=1$, $\phi_{n,\lambda}=\phi_{n,\lambda}(1)$ are called the degenerate Bell numbers, \par 
From \eqref{13}, we note that 
\begin{displaymath}
	\phi_{n,\lambda}(x)=\sum_{k=0}^{n}S_{2,\lambda}(n,k) x^{k},\quad (n\ge 0),\quad (\mathrm{see}\ [10,14]).
\end{displaymath}
Note that $\displaystyle \lim_{\lambda\rightarrow 0}\phi_{n,\lambda}(x)=\sum_{k=0}^{n}S_{2}(n,k)x^{k}\displaystyle$ are the ordinary Bell polynomials. \par 
For $k\in\mathbb{Z}$, the degenerate polylogarithm is defined by Kim-Kim as 
\begin{equation}
\mathrm{Li}_{k,\lambda}(t)=\sum_{n=1}^{\infty}\frac{\lambda^{n-1}(1)_{n,1/\lambda}}{(n-1)!n^{k}}(-1)^{n-1}t^{n},\quad (\mathrm{see}\ [8,13]).\label{14}
\end{equation}
Note that $\mathrm{Li}_{1,\lambda}(t)=-\log_{\lambda}(1-t)$. \par 
Carlitz considered the degenerate Bernoulli polynomials given by 
\begin{equation}
\frac{t}{e_{\lambda}(t)-1}e_{\lambda}^{x}(t)=\sum_{n=0}^{\infty}\beta_{n,\lambda}(x)\frac{t^{n}}{n!},\quad (\mathrm{see}\ [2,3]). \label{15}
\end{equation}
When $x=0$, $\beta_{n,\lambda}=\beta_{n,\lambda}(0)$ are called the degenerate Bernoulli numbers. \par 
He also defined the degenerate Euler polynomials given by 
\begin{equation}
\frac{2}{e_{\lambda}(t)+1}e_{\lambda}^{x}(t)=\sum_{n=0}^{\infty}\mathcal{E}_{n,\lambda}(x)\frac{t^{n}}{n!},\quad (\mathrm{see}\ [2]). \label{16}	
\end{equation}
Recently, Kim-Kim considered the degenerate Fubini polynomials defined by 
\begin{equation}
\frac{1}{1-x(e_{\lambda}(t)-1)}=\sum_{n=0}^{\infty}F_{n,\lambda}(x)\frac{t^{n}}{n!},\quad (\mathrm{see}\ [12]). \label{17}	
\end{equation}
By \eqref{11} and \eqref{17}, we easily get 
\begin{equation}
F_{n,\lambda}(x)=\sum_{k=0}^{n}S_{2,\lambda}(n,k)k!x^{k},\quad (n\ge 0),\quad (\mathrm{see}\ [12]).\label{18}
\end{equation}

\section{Some identities involving degenerate Stirling numbers associated with several degenerate polynomials and numbers}
In view of \eqref{6}, we consider the degenerate hyperharnomic numbers given by 
\begin{equation}
H_{n,\lambda}^{(1)}=H_{n,\lambda},\quad H_{n,\lambda}^{(r)}=\sum_{k=1}^{n}H_{k,\lambda}^{(r-1)},\quad (r\ge 2).\label{19}	
\end{equation}
By \eqref{5}, we easily get 
\begin{equation}
-\frac{\log_{\lambda}(1-t)}{1-t}=\sum_{n=1}^{\infty}H_{n,\lambda}t^{n},\quad (\mathrm{see}\ [9]).\label{20}
\end{equation}
In the view of \eqref{20}, we try to derive the generating function of the degenerate hyperharmonic numbers. \par 
From \eqref{20}, we note that 
\begin{align}
-\frac{\log_{\lambda}(1-t)}{(1-t)^{r}}&=\frac{1}{(1-t)^{r-1}}\bigg(-\frac{\log_{\lambda}(1-t)}{1-t}\bigg)=\frac{1}{(1-t)^{r-1}}\bigg(\sum_{n=1}^{\infty}H_{n,\lambda}t^{n}\bigg)\label{21} \\
&=\frac{1}{(1-t)^{r-2}}\sum_{n=1}^{\infty}\bigg(\sum_{k=1}^{n}H_{k,\lambda}\bigg)t^{n}=\frac{1}{(1-t)^{r-2}}\sum_{n=1}^{\infty}H_{n,\lambda}^{(2)}t^{n}\nonumber \\
&=\frac{1}{(1-t)^{r-3}}\sum_{n=1}^{\infty}\bigg(\sum_{k=1}^{n}H_{k,\lambda}^{(2)}\bigg)t^{n}=\frac{1}{(1-t)^{r-3}}\sum_{n=1}^{\infty}H_{n,\lambda}^{(3)}t^{n}.\nonumber
\end{align}
Continuing this process, we have 
\begin{equation}
-\frac{\log_{\lambda}(1-t)}{(1-t)^{r}}=\sum_{n=1}^{\infty}H_{n,\lambda}^{(r)}t^{n},\quad (r\ge 1),\quad H_{0,\lambda}^{(r)}=0. \label{22}
\end{equation}
Therefore, by \eqref{22}, we obtain the following theorem. 
\begin{theorem}
The generating function of the degenerate hyperharmonic numbers $H_{n,\lambda}^{(r)}$ is given by 
\begin{displaymath}
	-\frac{\log_{\lambda}(1-t)}{(1-t)^{r}}=\sum_{n=1}^{\infty}H_{n,\lambda}^{(r)}t^{n},\quad H_{0,\lambda}^{(r)}=0,\quad (r\ge 1). 
\end{displaymath}
\end{theorem}
Replacing $t$ by $1-e_{\lambda}(t)$ in \eqref{22}, we get 
\begin{align}
-te_{\lambda}^{-r}(t) &=\sum_{k=1}^{\infty}H_{k,\lambda}^{(r)}(-1)^{k}k!\frac{1}{k!}\big(e_{\lambda}(t)-1\big)^{k}\nonumber 	\\
&=\sum_{k=1}^{\infty}H_{k,\lambda}^{(r)}(-1)^{k}k!\sum_{n=k}^{\infty}S_{2,\lambda}(n,k)\frac{t^{n}}{n!} \label{23} \\
&=\sum_{n=1}^{\infty}\bigg(\sum_{k=1}^{n}H_{k,\lambda}^{(r)}(-1)^{k}k!S_{2,\lambda}(n,k)\bigg)\frac{t^{n}}{n!}.\nonumber
\end{align}
On the other hand, by \eqref{1}, we get 
\begin{align}
-te_{\lambda}^{-r}(t)&=-t\sum_{n=0}^{\infty}\frac{(-r)_{n,\lambda}}{n!}t^{n}=t\sum_{n=0}^{\infty}(-1)^{n-1}\frac{\langle r\rangle_{n,\lambda}}{n!}t^{n}\nonumber \\
&=\sum_{n=1}^{\infty}(-1)^{n}\langle r\rangle_{n-1,\lambda}\frac{t^{n}}{(n-1)!}=\sum_{n=1}^{\infty}(-1)^{n}\langle r\rangle_{n-1,\lambda}n	\frac{t^{n}}{n!}, \label{24}
\end{align}
where $\langle x\rangle_{0,\lambda}=1,\ \langle x\rangle_{n,\lambda}=x(x+\lambda)(x+2\lambda)\cdots(x+(n-1)\lambda),\ (n\ge 1)$. \par 
Therefore, by \eqref{23} and \eqref{24}, we obtain the following theorem. 
\begin{theorem}
For $n, r\in\mathbb{N}$ , we have 
\begin{displaymath}
	\sum_{k=1}^{n}(-1)^{k}H_{k,\lambda}^{(r)}k!S_{2,\lambda}(n,k)=(-1)^{n}\langle r\rangle_{n-1,\lambda}n. 
\end{displaymath}
In particular, for $r=1$, we get 
\begin{displaymath}
	\frac{1}{\langle 1\rangle_{n-1,\lambda}}\sum_{k=1}^{n}(-1)^{n-k}H_{k,\lambda}k!S_{2,\lambda}(n,k)=n.
\end{displaymath}
\end{theorem}
From \eqref{10}, we note that 
\begin{equation}
	\frac{1}{k!}\big(\log_{\lambda}(1+t)\big)^{k}=\sum_{n=k}^{\infty}S_{1,\lambda}(n,k)\frac{t^{n}}{n!},\quad (k\ge 0).\label{25}
\end{equation}
By \eqref{22} and \eqref{25}, we get 
\begin{align}	
\sum_{n=1}^{\infty}H_{n,\lambda}^{(r)}t^{n}&=-\frac{\log_{\lambda}(1-t)}{(1-t)^{r}}=-\log_{\lambda}(1-t)e_{\lambda}^{-r}\big(\log_{\lambda}(1-t)\big) \label{26}\\
&=-\log_{\lambda}(1-t)\sum_{k=0}^{\infty}(-1)^{k}\langle r\rangle_{k,\lambda}\frac{1}{k!}\big(\log_{\lambda}(1-t)\big)^{k} \nonumber \\
&=\sum_{k=1}^{\infty}(-1)^{k}\langle r\rangle_{k-1,\lambda}\frac{1}{(k-1)!}\big(\log_{\lambda}(1-t)\big)^{k}\nonumber \\
&=\sum_{k=1}^{\infty}(-1)^{k}\langle r\rangle_{k-1,\lambda}
k\sum_{n=k}^{\infty}S_{1,\lambda}(n,k)(-1)^{n}\frac{t^{n}}{n!} \nonumber \\
&=\sum_{n=1}^{\infty}\bigg(\sum_{k=1}^{n}(-1)^{n-k}\langle r\rangle_{k-1,\lambda}kS_{1,\lambda}(n,k)\bigg)\frac{t^{n}}{n!}. \nonumber
\end{align}
Therefore, by comparing the coefficients on both sides of \eqref{26}, we obtain the following theorem. 
\begin{theorem}
For $n, r\in\mathbb{N}$, we have 
\begin{displaymath}
H_{n,\lambda}^{(r)}=\frac{1}{n!}\sum_{k=1}^{n}(-1)^{n-k}\langle r\rangle_{k-1,\lambda}kS_{1,\lambda}(n,k).
\end{displaymath}
In particular, for $r=1$, we have
\begin{displaymath}
H_{n,\lambda}=\frac{1}{n!}\sum_{k=1}^{n}(-1)^{n-k}k\langle 1\rangle_{k-1,\lambda}S_{1,\lambda}(n,k).
\end{displaymath}
\end{theorem}
Replacing $t$ by $\log_{\lambda}(1-t)$ in \eqref{16}, we get 
\begin{align}
\frac{2}{2-t}(1-t)^{x}&=\sum_{k=0}^{\infty}\mathcal{E}_{k,\lambda}(x)\frac{1}{k!}\big(\log_{\lambda}(1-t)\big)^{k}\label{27}\\
&=\sum_{k=0}^{\infty}\mathcal{E}_{k,\lambda}(x)\sum_{n=k}^{\infty}S_{1,\lambda}(n,k)\frac{(-1)^{n}t^{n}}{n!} \nonumber \\
&=\sum_{n=0}^{\infty}\bigg(\sum_{k=0}^{n}(-1)^{n}S_{1,\lambda}(n,k)\mathcal{E}_{k,\lambda}(x)\bigg)\frac{t^{n}}{n!}.\nonumber	
\end{align}
By binomial expansion, we get 
\begin{align}
\frac{2}{2-t}(1-t)^{x}&=\frac{1}{1-\frac{t}{2}}(1-t)^{x}=\sum_{l=0}^{\infty}\bigg(\frac{1}{2}\bigg)^{l}t^{l}\sum_{k=0}^{\infty}\binom{x}{k}(-1)^{k}t^{k}\label{28}\\
&=\sum_{n=0}^{\infty}\bigg(\sum_{k=0}^{n}\binom{x}{k}\bigg(\frac{1}{2}\bigg)^{n-k}(-1)^{k}\bigg)t^{n}\nonumber.
\end{align}
Thus, by \eqref{27} and \eqref{28}, we get 
\begin{equation}
\sum_{k=0}^{n}S_{1,\lambda}(n,k)\mathcal{E}_{k,\lambda}(x)=n!\sum_{k=0}^{n}\binom{x}{k}\bigg(-\frac{1}{2}\bigg)^{n-k}.\label{29}
\end{equation}
Replacing $t$ by $1-e_{\lambda}(t)$ in \eqref{28}, we have 
\begin{align}
&\frac{2}{e_{\lambda}(t)+1}e_{\lambda}^{x}(t)=\sum_{k=0}^{\infty}\bigg(\sum_{j=0}^{k}\binom{x}{j}\bigg(\frac{1}{2}\bigg)^{k-j}(-1)^{j}\bigg)(1-e_{\lambda}(t))^{k}\label{30}	\\
&=\sum_{k=0}^{\infty}\bigg(\sum_{j=0}^{k}\binom{x}{j}\bigg(\frac{1}{2}\bigg)^{k-j}(-1)^{k-j}\bigg)k!\frac{1}{k!}\big(e_{\lambda}(t)-1\big)^{k} \nonumber \\
&=\sum_{k=0}^{\infty}\bigg(\sum_{j=0}^{k}\binom{x}{j}\bigg(\frac{1}{2}\bigg)^{k-j}(-1)^{k-j}k!\bigg)\sum_{n=k}^{\infty}S_{2,\lambda}(n,k)\frac{t^{n}}{n!} \nonumber \\
&=\sum_{n=0}^{\infty}\bigg(\sum_{k=0}^{n}\sum_{j=0}^{k}\binom{x}{j}\bigg(\frac{1}{2}\bigg)^{k-j}(-1)^{k-j}k!S_{2,\lambda}(n,k)\bigg)\frac{t^{n}}{n!}.\nonumber
\end{align}
By \eqref{16} and \eqref{30}, we get 
\begin{equation}
\mathcal{E}_{n,\lambda}(x)=\sum_{k=0}^{n}k!S_{2,\lambda}(n,k)\sum_{j=0}^{k}\binom{x}{j}\bigg(-\frac{1}{2}\bigg)^{k-j},\quad (n\ge 0). \label{31}
\end{equation}
In particular, for $x=\frac{1}{2}$, we obtain
\begin{align}
\mathcal{E}_{n,\lambda}\bigg(\frac{1}{2}\bigg)&=\sum_{k=0}^{n}k!S_{2,\lambda}(n,k)\sum_{j=0}^{k}\binom{\frac{1}{2}}{j}\bigg(-\frac{1}{2}\bigg)^{k-j} \label{32}\\
&=\sum_{k=0}^{n}k!S_{2,\lambda}(n,k)(-1)^{k}\sum_{j=0}^{k}\binom{2j}{j}\frac{1}{(1-2j)2^{k+j}}.\nonumber
\end{align}
Therefore, by \eqref{29} and \eqref{31}, we obtain the following theorem. 
\begin{theorem}
For $n\ge 0$, we have 
\begin{displaymath}
\sum_{k=0}^{n}S_{1,\lambda}(n,k)\mathcal{E}_{k,\lambda}(x)=n!\sum_{k=0}^{n}\binom{x}{k}\bigg(-\frac{1}{2}\bigg)^{n-k},
\end{displaymath}	
and 
\begin{displaymath}
\mathcal{E}_{n,\lambda}(x)=\sum_{k=0}^{n}k!S_{2,\lambda}(n,k)\sum_{j=0}^{k}\binom{x}{j}\bigg(-\frac{1}{2}\bigg)^{k-j}.
\end{displaymath}
\end{theorem}
From \eqref{3}, we note that 
\begin{align}
\frac{1}{t}\log_{\lambda}(1+t)&=\frac{1}{t}\sum_{n=1}^{\infty}\frac{\lambda^{n-1}(1)_{n,\frac{1}{\lambda}}}{n!}t^{n}\label{33} \\
&=\sum_{n=0}^{\infty}\frac{\lambda^{n}(1)_{n+1,\frac{1}{\lambda}}}{(n+1)!}t^{n}. \nonumber
\end{align}
Replacing $t$ by $e_{\lambda}(t)-1$ in \eqref{33}, we get 
\begin{align}
\sum_{n=0}^{\infty}\beta_{n,\lambda}\frac{t^{n}}{n!}&=\frac{t}{e_{\lambda}(t)-1}=\sum_{k=0}^{\infty}\frac{\lambda^{k}(1)_{k+1,1/\lambda}}{(k+1)!}k!\frac{1}{k!}(e_{\lambda}(t)-1)^{k} \label{34}\\
&=\sum_{k=0}^{\infty}\frac{\lambda^{k}(1)_{k+1,1/\lambda}}{k+1}\sum_{n=k}^{\infty}S_{2,\lambda}(n,k)\frac{t^{n}}{n!}\nonumber \\
&=\sum_{n=0}^{\infty}\sum_{k=0}^{n}\frac{\lambda^{k}(1)_{k+1,\frac{1}{\lambda}}}{k+1}S_{2,\lambda}(n,k)\frac{t^{n}}{n!}\nonumber \\
&=\sum_{n=0}^{\infty}\bigg(\sum_{k=0}^{n}\frac{\lambda^{k}(1)_{k+1,\frac{1}{\lambda}}}{k+1}S_{2,\lambda}(n,k)\bigg)\frac{t^{n}}{n!}.\nonumber
\end{align}
By comparing the coefficients on both sides of \eqref{34}, we get 
\begin{equation}
\beta_{n,\lambda}=\sum_{k=0}^{n}\frac{\lambda^{k}(1)_{k+1,\frac{1}{\lambda}}}{k+1}S_{2,\lambda}(n,k),\quad (n\ge 0). 	\label{35}
\end{equation}
From \eqref{15}, we note that 
\begin{align}
\beta_{n,\lambda}(x)&=\sum_{k=0}^{n}\binom{n}{k}\beta_{k,\lambda}(x)_{n-k,\lambda} \label{36}\\
&=\sum_{k=0}^{n}\binom{n}{k}(x)_{k,\lambda}\beta_{n-k,\lambda},\quad (n\ge 0). \nonumber
\end{align}
By \eqref{35} and \eqref{36}, we see that 
\begin{align}
\beta_{n,\lambda}(x)&=\sum_{k=0}^{n}\binom{n}{k}(x)_{n-k,\lambda}\beta_{k,\lambda}\label{37}\\
 &=\sum_{k=0}^{n}\binom{n}{k}(x)_{n-k,\lambda}\sum_{j=0}^{k}\frac{\lambda^{j}(1)_{j+1,\frac{1}{\lambda}}}{j+1}S_{2,\lambda}(k,j)\nonumber\\
&=\sum_{j=0}^{n}\frac{\lambda^{j}(1)_{j+1,\frac{1}{\lambda}}}{j+1}\sum_{k=j}^{n}\binom{n}{k}S_{2,\lambda}(k,j)(x)_{n-k,\lambda}\nonumber.
\end{align}

Therefore, by \eqref{37}, we obtain the following theorem. 
\begin{theorem}
For $n\ge 0$, we have 
\begin{displaymath}
	\beta_{n,\lambda}(x)=\sum_{j=0}^{n}\frac{\lambda^{j}(1)_{j+1,\frac{1}{\lambda}}}{j+1}\sum_{k=j}^{n}\binom{n}{k}S_{2,\lambda}(k,j)(x)_{n-k,\lambda}.
\end{displaymath}	
\end{theorem}
Replacing $t$ by $\log_{\lambda}(1+t)$ in \eqref{15}, we have 
\begin{align}
\frac{\log_{\lambda}(1+t)}{t}(1+t)^{x}&=\sum_{k=0}^{\infty}\beta_{k,\lambda}(x)\frac{1}{k!}\big(\log_{\lambda}(1+t)\big)^{k} \label{38}	\\
&=\sum_{k=0}^{\infty}\beta_{k,\lambda}(x)\sum_{n=k}^{\infty}S_{1,\lambda}(n,k)\frac{t^{n}}{n!} \nonumber\\
&=\sum_{n=0}^{\infty}\bigg(\sum_{k=0}^{n}\beta_{k,\lambda}(x)S_{1,\lambda}(n,k)\bigg)\frac{t^{n}}{n!}. \nonumber
\end{align}
On the other hand, by \eqref{3}, we get 
\begin{align}
\frac{\log_{\lambda}(1+t)}{t}(1+t)^{x}&=\frac{1}{t}\sum_{n=1}^{\infty}\frac{\lambda^{n-1}(1)_{n,\frac{1}{\lambda}}}{n!}t^{n} (1+t)^x\label{39}\\
&=\sum_{l=0}^{\infty}\frac{\lambda^{l}(1)_{l+1,\frac{1}{\lambda}}}{(l+1)!}t^{l}\sum_{k=0}^{\infty}\binom{x}{k}t^{k} \nonumber \\
&=\sum_{n=0}^{\infty}\bigg(n!\sum_{k=0}^{n}\binom{x}{k}\frac{\lambda^{n-k}(1)_{n-k+1,\frac{1}{\lambda}}}{(n-k+1)!}\bigg)\frac{t^{n}}{n!}. \nonumber
\end{align}
Therefore, by \eqref{38} and \eqref{39}, we obtain the following theorem. 
\begin{theorem}
For $n\ge 0$, we have 
\begin{displaymath}
	\sum_{k=0}^{n}\beta_{k,\lambda}(x)S_{1,\lambda}(n,k)=n!\sum_{k=0}^{n}\binom{x}{k}\frac{\lambda^{n-k}(1)_{n-k+1,\frac{1}{\lambda}}}{(n-k+1)!}.
\end{displaymath}
\end{theorem}
In \eqref{39}, by replacing $t$ by $e_{\lambda}(t)-1$, we get 
\begin{align}
&\sum_{n=0}^{\infty}\beta_{n,\lambda}(x)\frac{t^{n}}{n!}=\frac{t}{e_{\lambda}(t)-1}e_{\lambda}^{x}(t)\label{40}\\
&=\sum_{k=0}^{\infty}k!	\sum_{j=0}^{k}\binom{x}{j}\frac{\lambda^{k-j}}{(k-j+1)!}(1)_{k-j+1,\frac{1}{\lambda}}\cdot\frac{1}{k!}\big(e_{\lambda}(t)-1\big)^{k}\nonumber \\
&=\sum_{k=0}^{\infty}k!\sum_{j=0}^{k}\binom{x}{j}\frac{\lambda^{k-j}}{(k-j+1)!}(1)_{k-j+1,\frac{1}{\lambda}}\sum_{n=k}^{\infty}S_{2,\lambda}(n,k)\frac{t^{n}}{n!} \nonumber \\
&=\sum_{n=0}^{\infty}\bigg(\sum_{k=0}^{n}k!S_{2,\lambda}(n,k)\sum_{j=0}^{k}\frac{\lambda^{k-j}}{(k-j+1)!}(1)_{k-j+1,\frac{1}{\lambda}}\binom{x}{j}\bigg)\frac{t^{n}}{n!}.\nonumber
\end{align}
Therefore, by comparing the coefficients on both sides of \eqref{40}, we obtain the following theorem. 
\begin{theorem}
For $n\ge 0$, we have 
\begin{displaymath}
\beta_{n,\lambda}(x)=\sum_{k=0}^{n}k!S_{2,\lambda}(n,k)\sum_{j=0}^{k}\binom{x}{j}\frac{\lambda^{k-j}}{(k-j+1)!}(1)_{k-j+1,\frac{1}{\lambda}}.
\end{displaymath}	
In particular, for $x=0$, we have
\begin{displaymath}
\beta_{n,\lambda}=\sum_{k=0}^{n}k!S_{2,\lambda}(n,k)\frac{\lambda^{k}}{(k+1)!}(1)_{k+1,\frac{1}{\lambda}}.
\end{displaymath}
\end{theorem}
Let us take $k=2$ in \eqref{14}. Then 
\begin{equation}
\mathrm{Li}_{2,\lambda}(t)=\sum_{n=1}^{\infty}\frac{\lambda^{n-1}(1)_{n,\frac{1}{\lambda}}}{(n-1)!n^{2}}(-1)^{n-1}t^{n},\quad (\mathrm{see}\ [13]).\label{41}
\end{equation}
Note that 
\begin{align}
-\int_{0}^{t}\frac{1}{x}\log_{\lambda}(1-x)dx&=\int_{0}^{t}\sum_{n=1}^{\infty}\frac{\lambda^{n-1}}{n!}(1)_{n,\frac{1}{\lambda}}(-1)^{n-1}x^{n-1}dx \label{42} \\
&=\sum_{n=1}^{\infty}\frac{\lambda^{n-1}}{(n-1)!n^{2}}(1)_{n,\frac{1}{\lambda}}(-1)^{n-1}t^{n}=\mathrm{Li}_{2,\lambda}(t).\nonumber	
\end{align}
Thus, by \eqref{42}, we get 
\begin{align}
\mathrm{Li}_{2}(1-e_{\lambda}(-t))&=-\int_{0}^{1-e_{\lambda}(-t)}\frac{1}{x}\log_{\lambda}(1-x)dx \label{43}\\
&=\int_{0}^{t}\frac{-t}{e_{\lambda}(-t)-1}e_{\lambda}^{1-\lambda}(-t)dt\nonumber \\
&=\sum_{n=1}^{\infty}(-1)^{n-1}\beta_{n-1,\lambda}(1-\lambda)\frac{t^{n}}{n!}.\nonumber	
\end{align}
From \eqref{42}, we note that 
\begin{align}
&\mathrm{Li}_{2,\lambda}\bigg(-\frac{t}{1-t}\bigg)=-\int_{0}^{\frac{-t}{1-t}}\frac{1}{x}\log_{\lambda}(1-x)dx\label{44}	\\
&=-\int_{0}^{t}\bigg(-\frac{1-t}{t}\bigg)\bigg(\log_{\lambda}\frac{1}{1-t}\bigg)\frac{-1}{(1-t)^{2}}dt\nonumber\\
&=-\int_{0}^{t}\frac{1}{t} \bigg(-\frac{\log_{-\lambda}(1-t)}{1-t}\bigg)dt=-\int_{0}^{t}\frac{1}{t}\sum_{n=1}^{\infty}H_{n,-\lambda}t^{n}dt \nonumber \\
&=-\sum_{n=1}^{\infty}(n-1)!H_{n,-\lambda}\frac{t^{n}}{n!}.\nonumber
\end{align}
Replacing $t$ by $-\log_{\lambda}\frac{1}{1-t}$ in \eqref{43}, we have 
\begin{align}
\mathrm{Li}_{2,\lambda}\bigg(\frac{-t}{1-t}\bigg)&=\sum_{k=1}^{\infty}(-1)^{k-1}\beta_{k-1,\lambda}(1-\lambda)\frac{1}{k!}\big(-\log_{\lambda}\big(\frac{1}{1-t}\big)\big)^{k}\label{45} \\
&=\sum_{k=1}^{\infty}(-1)^{k-1}\beta_{k-1,\lambda}(1-\lambda)\frac{1}{k!}\big(\log_{-\lambda}(1-t)\big)^{k}\nonumber \\
&=\sum_{k=1}^{\infty}(-1)^{k-1}\beta_{k-1,\lambda}(1-\lambda)\sum_{n=k}^{\infty}S_{1,-\lambda}(n,k)\frac{(-t)^{n}}{n!} \nonumber \\
&=\sum_{n=1}^{\infty}\bigg(\sum_{k=1}^{n}(-1)^{n-k-1}\beta_{k-1,\lambda}(1-\lambda)S_{1,-\lambda}(n,k)\bigg)\frac{t^{n}}{n!} .\nonumber
\end{align}
Therefore, by \eqref{44} and \eqref{45}, we obtain the following theorem. 
\begin{theorem}
	For $n\in\mathbb{N}$, we have 
	\begin{displaymath}
		(n-1)!H_{n,-\lambda}=\sum_{k=1}^{n}(-1)^{n-k}\beta_{k-1,\lambda}(1-\lambda)S_{1,-\lambda}(n,k).
	\end{displaymath}
\end{theorem}
In \eqref{44}, by replacing $t$ by $1-e_{\lambda}(t)$, we get 
\begin{align}
\mathrm{Li}_{2,\lambda}\bigg(\frac{e_{\lambda}(t)-1)}{e_{\lambda}(t)}\bigg)&=-\sum_{k=1}^{\infty}(k-1)!H_{k,-\lambda}\frac{1}{k!}\big(1-e_{\lambda}(t)\big)^{k}\label{46}\\
&=-\sum_{k=1}^{\infty}(k-1)!H_{k,-\lambda}(-1)^{k}\sum_{n=k}^{\infty}S_{2,\lambda}(n,k)\frac{t^{n}}{n!}\nonumber \\
&=\sum_{n=1}^{\infty}\bigg(\sum_{k=1}^{n}(k-1)!(-1)^{k-1}H_{k,-\lambda}S_{2,\lambda}(n,k)\bigg)\frac{t^{n}}{n!}.\nonumber	
\end{align}
On the other hand, by \eqref{42}, we get 
\begin{align}
&\mathrm{Li}_{2,\lambda}\bigg(\frac{e_{\lambda}(t)-1}{e_{\lambda}(t)}\bigg)=\mathrm{Li}_{2,\lambda}\bigg(1-\frac{1}{e_{\lambda}(t)}\bigg)=\mathrm{Li}_{2,\lambda}\big(1-e_{-\lambda}(-t)\big) \label{47}\\
&=-\int_{0}^{1-e_{-\lambda}(-t)}\frac{1}{x}\log_{\lambda}(1-x)dx 
=\int_{0}^{t}\frac{-t}{e_{-\lambda}(-t)-1}e_{-\lambda}^{2\lambda+1}(-t)dt\nonumber \\
&=\int_{0}^{t}\sum_{n=0}^{\infty}\beta_{n,-\lambda}(2\lambda+1)(-1)^{n}\frac{t^{n}}{n!}dt =\sum_{n=1}^{\infty}\beta_{n-1,-\lambda}(2\lambda+1)(-1)^{n-1}\frac{t^{n}}{n!}. \nonumber
\end{align}
Therefore, by \eqref{46} and \eqref{47}, we obtain the following theorem. 
\begin{theorem}
For $n\in\mathbb{N}$, we have 
\begin{displaymath}
\beta_{n-1,-\lambda}(2\lambda+1)=\sum_{k=1}^{n}(k-1)!(-1)^{n-k}H_{k,-\lambda}S_{2,\lambda}(n,k).
\end{displaymath}
\end{theorem}
Replacing $t$ by $-\log_{\lambda}(1-t)$ in \eqref{43}, we have 
\begin{align}
&\sum_{n=1}^{\infty}\frac{\lambda^{n-1}(1)_{n,\frac{1}{\lambda}}}{(n-1)!n^{2}}(-1)^{n-1}t^{n}=\mathrm{Li}_{2,\lambda}(t)=\sum_{k=1}^{\infty}(-1)^{k-1}\beta_{k-1,\lambda}(1-\lambda)\frac{(-\log_{\lambda}(1-t))^{k}}{k!}\label{48} \\
&=-\sum_{k=1}^{\infty}\beta_{k-1,\lambda}(1-\lambda)\sum_{n=k}^{\infty}S_{1,\lambda}(n,k)(-1)^{n}\frac{t^{n}}{n!}=\sum_{n=1}^{\infty}\bigg(\sum_{k=1}^{n}(-1)^{n-1}\beta_{k-1,\lambda}(1-\lambda)S_{1,\lambda}(n,k)\bigg)\frac{t^{n}}{n!},\nonumber
\end{align}
Therefore, by comparing the coefficients on both sides of \eqref{48}, we obtain the following theorem. 
\begin{theorem}
For $n\in\mathbb{N}$, we have 
\begin{displaymath}
\frac{\lambda^{n-1}}{n}(1)_{n,\frac{1}{\lambda}}=\sum_{k=1}^{n}\beta_{k-1,\lambda}(1-\lambda)S_{1,\lambda}(n,k).
\end{displaymath}
\end{theorem}
From \eqref{13}, we note that 
\begin{align}
\phi_{n,\lambda}(x)&=\frac{1}{e^{x}}\sum_{k=0}^{\infty}\frac{(k)_{n,\lambda}}{k!}x^{k} 
=\frac{1}{e^{x}}\bigg(x\frac{d}{dx}\bigg)_{n,\lambda}e^{x},\quad (n\ge 1).\label{49}	
\end{align}
Thus, by \eqref{49}, we get 
\begin{equation}
\phi_{n,\lambda}(x)=e^{-x}\bigg(x\frac{d}{dx}\bigg)_{n,\lambda}e^{x}=\sum_{k=1}^{n}S_{2,\lambda}(n,k)x^{k},\quad (n\ge 1).\label{50}	
\end{equation}
From \eqref{13}, we have 
\begin{equation}
S_{2,\lambda}(n+1,k)=S_{2,\lambda}(n,k-1)+(k-n\lambda)S_{2,\lambda}(n,k),\label{51}
\end{equation}
where $n,k$ are nonnegative integers with $n \ge k$. \par 
By \eqref{50} and \eqref{51}, we get 
\begin{align}
\phi_{n+1,\lambda}(x)&=\sum_{k=0}^{n+1}S_{2,\lambda}(n+1,k)x^{k}=\sum_{k=0}^{n+1}\big\{S_{2,\lambda}(n,k-1)+(k-n\lambda)S_{2,\lambda}(n,k)\big\}x^{k}\label{52}\\
&=\sum_{k=1}^{n+1}S_{2,\lambda}(n,k-1)x^{k}+\sum_{k=0}^{n+1}S_{2,\lambda}(n,k)(k-n\lambda)x^{k}\nonumber \\
&=x\sum_{k=0}^{n}S_{2,\lambda}(n,k)x^{k}+\Big(x\frac{d}{dx}-n\lambda\Big)\sum_{k=0}^{n}S_{2,\lambda}(n,k)x^{k}\nonumber \\
&=x\phi_{n,\lambda}(x)+\Big(x\frac{d}{dx}-n\lambda\Big)\phi_{n,\lambda}(x).\nonumber
\end{align}
By \eqref{13}, we easily get 
\begin{align}
\sum_{n=0}^{\infty}\phi_{n+1,\lambda}(x)\frac{t^{n}}{n!}&=\frac{d}{dt}\sum_{n=0}^{\infty}\phi_{n,\lambda}(x)\frac{t^{n}}{n!} \label{53} \\
&=\frac{d}{dt}e^{x(e_{\lambda}(t)-1)}=xe_{\lambda}^{1-\lambda}(t)e^{x(e_{\lambda}(t)-1)}\nonumber \\
&=\sum_{n=0}^{\infty}\bigg(x\sum_{k=0}^{n}\binom{n}{k}(1-\lambda)_{n-k,\lambda}\phi_{k,\lambda}(x)\bigg)\frac{t^{n}}{n!}. \nonumber
\end{align}
Therefore, by \eqref{52} and \eqref{53}, we obtain the following theorem. 
\begin{theorem}
For $n\in\mathbb{N}$, we have 
\begin{align*}
\phi_{n+1,\lambda}(x)&=x\phi_{n,\lambda}(x)+\bigg(x\frac{d}{dx}-n\lambda\bigg)\phi_{n,\lambda}(x)\\
&=x\sum_{k=0}^{n}\binom{n}{k}(1-\lambda)_{n-k,\lambda}\phi_{k,\lambda}(x).
\end{align*}
\end{theorem}
For $p \ge 0$, we observe that 
\begin{align}
&\sum_{k=1}^{p+1}S_{2,\lambda}(p+1,k)x^{k}=\phi_{p+1,\lambda}	(x)=e^{-x}\sum_{n=1}^{\infty}(n)_{p+_1,\lambda}\frac{x^{n}}{n!} \label{54}\\
&=e^{-x}\sum_{n=0}^{\infty}(n+1)_{p+1,\lambda}\frac{x^{n+1}}{(n+1)!}=e^{-x}\sum_{n=0}^{\infty}(n+1-\lambda)_{p,\lambda}\frac{x^{n+1}}{n!}. \nonumber
\end{align}
Thus, by \eqref{54}, we get 
\begin{align}
&\sum_{k=1}^{p+1}S_{2,\lambda}(p+1,k)x^{k-1}=e^{-x}\sum_{n=0}^{\infty}(n+1-\lambda)_{p,\lambda}\frac{x^{n}}{n!}\label{55}	\\
&=e^{-x}\sum_{n=0}^{\infty}\Big\{(1-\lambda)_{p,\lambda}+(2-\lambda)_{p,\lambda}+\cdots+(n+1-\lambda)_{p,\lambda}-((1-\lambda)_{p,\lambda}+\cdots+(n-\lambda)_{p,\lambda}\Big\}\frac{x^{n}}{n!}\nonumber\\
&=e^{-x}\sum_{n=1}^{\infty}\big((1-\lambda)_{p,\lambda}+(2-\lambda)_{p,\lambda}+\cdots+(n-\lambda)_{p,\lambda}\big)\frac{x^{n-1}}{(n-1)!} \nonumber \\
&\quad -e^{-x}\sum_{n=0}^{\infty}\big((1-\lambda)_{p,\lambda}+(2-\lambda)_{p,\lambda}+\cdots+(n-\lambda)_{p,\lambda}\big)\frac{x^{n}}{n!} \nonumber\\
&=\frac{d}{dx}\bigg(e^{-x}\sum_{n=0}^{\infty}\big((1-\lambda)_{p,\lambda}+(2-\lambda)_{p,\lambda}+\cdots+(n-\lambda)_{p,\lambda}\big)\frac{x^{n}}{n!}\bigg)\nonumber.
\end{align}
From \eqref{55}, we note that 
\begin{equation}
\frac{d}{dx}\sum_{k=1}^{p+1}S_{2,\lambda}(p+1,k)\frac{x^{k}}{k}=\frac{d}{dx}\bigg(e^{-x}\sum_{n=0}^{\infty}\big((1-\lambda)_{p,\lambda}+(2-\lambda)_{p,\lambda}+\cdots+(n-\lambda)_{p,\lambda}\big)\frac{x^{n}}{n!}\bigg).\label{56}
\end{equation}
Thus, by \eqref{56}, we get 
\begin{equation}
\sum_{k=1}^{p+1}S_{2,\lambda}(p+1,k)\frac{x^{k}}{k}=e^{-x}\sum_{n=0}^{\infty}\big((1-\lambda)_{p,\lambda}+(2-\lambda)_{p,\lambda}+\cdots+(n-\lambda)_{p,\lambda}\big)\frac{x^{n}}{n!}.\label{57}	
\end{equation}
Now, we observe that 
\begin{align}
&\sum_{p=0}^{\infty}\sum_{k=0}^{n-1}(k-\lambda+1)_{p,\lambda}\frac{t^{p}}{p!}=\sum_{k=0}^{n-1}e_{\lambda}^{k+1-\lambda}(t) \label{58} \\
&=e_{\lambda}^{1-\lambda}(t)\frac{e_{\lambda}^{n}(t)-1}{e_{\lambda}(t)-1}=\sum_{p=0}^{\infty}\bigg(\frac{\beta_{p+1}(n+1-\lambda)-\beta_{p+1,\lambda}(1-\lambda)}{p+1}\bigg)\frac{t^{p}}{p!}.\nonumber
\end{align}
Comparing the coefficients on both sides of \eqref{58}, we have 
\begin{align}
\sum_{k=1}^{n}(k-\lambda)_{p,\lambda}&=\frac{1}{p+1}\big(\beta_{p+1,\lambda}(n+1-\lambda)-\beta_{p+1,\lambda}(1-\lambda)\big) \label{59}\\
&=\frac{1}{p+1}\bigg\{\sum_{j=0}^{p+1}\binom{p+1}{j}\beta_{j,\lambda}(1-\lambda)(n)_{p+1-j,\lambda}-\beta_{p+1,\lambda}(1-\lambda)\bigg\}\nonumber\\
&=\frac{1}{p+1}\sum_{j=0}^{p}\binom{p+1}{j}\beta_{j,\lambda}(1-\lambda)(n)_{p+1-j,\lambda}.\nonumber
\end{align}
So, by \eqref{57} and \eqref{59}, we get 
\begin{align}
\sum_{k=1}^{p+1}S_{2,\lambda}(p+1,k)\frac{x^{k}}{k}&=e^{-x}\sum_{n=0}^{\infty}\big((1-\lambda)_{p,\lambda}+\cdots+(n-\lambda)_{p,\lambda}\big)\frac{x^{n}}{n!}.\label{60}\\
&=e^{-x}\sum_{n=0}^{\infty}\bigg(\frac{1}{p+1}\sum_{j=0}^{p}\binom{p+1}{j}\beta_{j,\lambda}(1-\lambda)(n)_{p+1-j,\lambda}\bigg)\frac{x^{n}}{n!} 	\nonumber \\
&=\frac{1}{p+1}\sum_{j=0}^{p}\binom{p+1}{j}\beta_{j,\lambda}(1-\lambda)\bigg(e^{-x}\sum_{n=0}^{\infty}(n)_{p+1-j,\lambda}\frac{x^{n}}{n!}\bigg) \nonumber \\
&=\frac{1}{p+1}\sum_{j=0}^{p}\binom{p+1}{j}\beta_{j,\lambda}(1-\lambda)\phi_{p+1-j,\lambda}(x),\quad (p\ge 0).\nonumber
\end{align}
Therefore, by \eqref{60}, we obtain the following theorem. 
\begin{theorem}
For $n\in\mathbb{Z}$ with $n\ge 1$, we have 
\begin{displaymath}
\sum_{k=1}^{n}S_{2,\lambda}(n,k)\frac{x^{k}}{k}=\frac{1}{n}\sum_{j=0}^{n-1}\binom{n}{j}\beta_{j,\lambda}(1-\lambda)\phi_{n-j,\lambda}(x).
\end{displaymath}
In particular, for $x=1$, we get 
\begin{displaymath}
\sum_{k=1}^{n}\frac{S_{2,\lambda}(n,k)}{k}=\frac{1}{n!}\sum_{j=0}^{n-1}\binom{n}{j}\beta_{j,\lambda}(1-\lambda)\phi_{n-j,\lambda}.
\end{displaymath}
\end{theorem}
From Theorem 12, we note that 
\begin{align}
\sum_{k=1}^{n}&S_{2,\lambda}(n,k)\frac{x^{k}}{k^{2}}=\frac{1}{n}\sum_{j=0}^{n-1}\binom{n}{j}\beta_{j,\lambda}(1-\lambda)	\int_{0}^{x}\frac{\phi_{n-j,\lambda}(t)}{t}dt \label{61} \\
&= \frac{1}{n}\sum_{j=0}^{n-1}\binom{n}{j}\beta_{j,\lambda}(1-\lambda)	\sum_{l=1}^{n-j}S_{2,\lambda}(n-j,l)\frac{x^{l}}{l} \nonumber \\
&= \frac{1}{n}\sum_{j=0}^{n-1}\binom{n}{j}\beta_{j,\lambda}(1-\lambda)\frac{1}{n-j}\sum_{m=1}^{n-j}\binom{n-j}{m}\beta_{n-j-m,\lambda}(1-\lambda)\phi_{m,\lambda}(x).\nonumber
\end{align}
Thus, by \eqref{61}, we get 
\begin{displaymath}
	\sum_{k=1}^{n}S_{2,\lambda}(n,k)\frac{x^{k}}{k^{2}} =\frac{1}{n}\sum_{j=0}^{n-1}\binom{n}{j}\beta_{j,\lambda}(1-\lambda)\frac{1}{n-j}\sum_{m=1}^{n-j}\binom{n-j}{m}\beta_{n-j-m,\lambda}(1-\lambda)\phi_{m,\lambda}(x),
\end{displaymath}
for $n\ge 1$. \par 
From \eqref{17} and \eqref{18}, we have 
\begin{equation}
\int_{0}^{x}\frac{F_{n,\lambda}(t)}{t}dt=\sum_{k=1}^{n}S_{2,\lambda}(n,k)(k-1)!x^{k},\quad (n\ge 1).\label{62}
\end{equation}
By \eqref{51}, we get 
\begin{align*}
	&\sum_{k=1}^{n+1}S_{2,\lambda}(n+1,k)(k-1)!x^{k}=\sum_{k=1}^{n+1}\big\{S_{2,\lambda}(n,k-1)+(k-n\lambda)S_{2,\lambda}(n,k)\big\}(k-1)!x^{k}\\
	&=\sum_{k=1}^{n}k!S_{2,\lambda}(n,k)x^{k}+x\sum_{k=1}^{n}S_{2,\lambda}(n,k)k!x^{k}-n\lambda \sum_{k=1}^{n}S_{2,\lambda}(n,k)(k-1)!x^{k} \\
	&=(1+x)F_{n,\lambda}(x)-n\lambda\int_{0}^{x}\frac{F_{n,\lambda}(t)}{t}dt,\quad (n\ge 1). 
\end{align*}
Therefore, by \eqref{62}, we obtain the following theorem.
\begin{theorem}
	For $n\in\mathbb{N}$, we have 
	\begin{displaymath}
		(1+x)F_{n,\lambda}(x)=\sum_{k=1}^{n+1}S_{2,\lambda}(n+1,k)(k-1)!x^{k}+n\lambda\sum_{k=1}^{n}S_{2,\lambda}(n,k)(k-1)!x^{k}.
	\end{displaymath}
	Moreover, we have
	\begin{displaymath}
		\int_{0}^{x}\frac{F_{n+1,\lambda}(t)}{t}dt+n\lambda\int_{0}^{x}\frac{F_{n,\lambda}(t)}{t}dt=(1+x)F_{n,\lambda}(x).
	\end{displaymath}
\end{theorem}
From Theorem 13, we have 
\begin{align}
\sum_{k=1}^{n+1}S_{2,\lambda}(n+1,k)k!x^{k-1}&=\frac{d}{dx}\sum_{k=1}^{n+1}S_{2,\lambda}(n+1,k)(k-1)!x^{k}\label{63}\\
&=-\frac{n\lambda}{x}\sum_{k=1}^{n}S_{2,\lambda}(n,k)k!x^{k}+F_{n,\lambda}(x)+(1+x)F_{n,\lambda}^{\prime}(x),\nonumber
\end{align}
where $F_{n,\lambda}^{\prime}(x)=\frac{d}{dx}F_{n,\lambda}(x)$. \par 
Thus, by \eqref{63}, we get 
\begin{equation}
\sum_{k=1}^{n+1}S_{2,\lambda}(n+1,k)k!x^{k}=-n\lambda	F_{n,\lambda}(x)+xF_{n,\lambda}(x)+(x+x^{2})F_{n,\lambda}^{\prime}(x). \label{64}
\end{equation}
Therefore, by \eqref{18} and \eqref{64}, we obtain the following theorem. 
\begin{theorem}
	For $n\in\mathbb{N}$, we have 
	\begin{displaymath}
		F_{n+1,\lambda}(x)=(x-n\lambda)F_{n,\lambda}+(x+x^{2})F_{n,\lambda}^{\prime}(x).
	\end{displaymath}
\end{theorem}
\emph{Remark.} By Theorem 13, we easily get 
\begin{displaymath}
	\sum_{k=1}^{n+1}S_{2,\lambda}(n+1,k)(k-1)!+n\lambda\sum_{k=1}^{n}S_{2,\lambda}(n,k)(k-1)!=\left\{\begin{array}{ccc}
		2F_{n,\lambda}(1), & \textrm{ if $n\ge 1$,}\\
		1, & \textrm{if $n=0$.}
	\end{array}\right.
\end{displaymath}
We observe that 
\begin{align}
&\sum_{n=0}^{\infty}\bigg(\sum_{k=0}^{\infty}x^{k}(k)_{n,\lambda}\bigg)\frac{t^{n}}{n!}=\sum_{k=0}^{\infty}x^{k}e_{\lambda}^{k}(t)=\frac{1}{1-xe_{\lambda}(t)} \label{65}\\
&=\frac{1}{1-x-x(e_{\lambda}(t)-1)}=\frac{1}{1-x}\frac{1}{1-\frac{x}{1-x}(e_{\lambda}(t)-1)}=\frac{1}{1-x}\sum_{n=0}^{\infty}F_{n,\lambda}\bigg(\frac{x}{1-x}\bigg)\frac{t^{n}}{n!}\nonumber.
\end{align}
Comparing the coefficients on both sides of \eqref{65}, we have 
\begin{equation}
\frac{1}{1-x}F_{n,\lambda}\bigg(\frac{x}{1-x}\bigg)=\sum_{k=0}^{\infty}x^{k}(k)_{n,\lambda},\quad (n\ge 0).\label{66}	
\end{equation}
Let us take $x=\frac{1}{2}$ in \eqref{66}. Then we get 
\begin{equation}
F_{n,\lambda}(1)=\sum_{k=0}^{\infty}(k)_{n,\lambda}\bigg(\frac{1}{2}\bigg)^{k+1},\quad (n\ge 0).\label{67}
\end{equation}
Let 
\begin{equation}
f(t)=(1+t)\log_{\lambda}(1+t)-t.\label{68}
\end{equation}
Then, by \eqref{3}, we get 
\begin{align}
f(t) &=t\log_{\lambda}(1+t)+\log_{\lambda}(1+t)-t \label{69}\\
&= t\sum_{n=1}^{\infty}\frac{\lambda^{n-1}(1)_{n,\frac{1}{\lambda}}}{n!}t^{n}+\sum_{n=2}^{\infty}\frac{\lambda^{n-1}(1)_{n,\frac{1}{\lambda}}}{n!}t^{n} \nonumber \\
&=t^{2}\sum_{n=0}^{\infty}\frac{\lambda^{n}(1)_{n+1,\frac{1}{\lambda}}}{n+1}\frac{t^{n}}{n!}+t^{2}\sum_{n=0}^{\infty}\frac{\lambda^{n+1}(1)_{n+1,\frac{1}{\lambda}}}{(n+2)!}\bigg(1-(n+1)\frac{1}{\lambda}\bigg)t^{n}\nonumber \\
&=t^{2}\sum_{n=0}^{\infty}\frac{\lambda^{n}(1)_{n+1,\frac{1}{\lambda}}}{n+1}\frac{t^{n}}{n!}-t^{2}\sum_{n=0}^{\infty}\frac{\lambda^{n}(1)_{n+1,\frac{1}{\lambda}}}{n+2}\frac{t^{n}}{n!}+t^{2}\sum_{n=0}^{\infty}\frac{\lambda^{n+1}(1)_{n+1,\frac{1}{\lambda}}}{(n+2)!}t^{n}\nonumber \\
&=t^{2}\sum_{n=0}^{\infty}\frac{\lambda^{n}(1)_{n+1,\frac{1}{\lambda}}}{(n+1)(n+2)}\frac{t^{n}}{n!}+\lambda t^{2}\sum_{n=0}^{\infty}\frac{\lambda^{n}(1)_{n+1,\frac{1}{\lambda}}}{(n+2)!}t^{n} \nonumber \\
&=\sum_{n=2}^{\infty}\lambda^{n-2}(1)_{n-1,\frac{1}{\lambda}}\frac{t^{n}}{n!}+\lambda\sum_{n=2}^{\infty}\lambda^{n-2}(1)_{n-1,\frac{1}{\lambda}}\frac{t^{n}}{n!} \nonumber \\
&=\sum_{n=2}^{\infty}(1+\lambda)\lambda^{n-2}(1)_{n-1,\frac{1}{\lambda}}\frac{t^{n}}{n!}. \nonumber
\end{align}
In \eqref{69}, by replacing $t$ by $e_{\lambda}(t)-1$, we get 
\begin{align}
&te_{\lambda}(t)-e_{\lambda}(t)+1=f(e_{\lambda}(t)-1)=(1+\lambda)\sum_{k=2}^{\infty}(1)_{k-1,\frac{1}{\lambda}}\lambda^{k-2}\frac{1}{k!}\big(e_{\lambda}(t)-1\big)^{k}\label{70}\\
&=(1+\lambda)\sum_{k=2}^{\infty}(1)_{k-1,\frac{1}{\lambda}}\lambda^{k-2}\sum_{n=k}^{\infty}S_{2,\lambda}(n,k)\frac{t^{n}}{n!} \nonumber \\
&=\sum_{n=2}^{\infty}\bigg((1+\lambda)\sum_{k=2}^{n}(1)_{k-1,\frac{1}{\lambda}}\lambda^{k-2}S_{2,\lambda}(n,k)\bigg)\frac{t^{n}}{n!}. \nonumber
\end{align}
On the other hand, by \eqref{1}, we get 
\begin{align}
te_{\lambda}(t)-e_{\lambda}(t)+1&=t\sum_{n=0}^{\infty}\frac{(1)_{n,\lambda}}{n!}t^{n}-\sum_{n=0}^{\infty}\frac{(1)_{n,\lambda}}{n!}t^{n}+1\label{71} \\
&=\sum_{n=0}^{\infty}\frac{(1)_{n,\lambda}}{n!}t^{n+1}-\sum_{n=1}^{\infty}\frac{(1)_{n,\lambda}}{n!}t^{n}\nonumber\\
&=\sum_{n=1}^{\infty}(1)_{n-1,\lambda}\frac{t^{n}}{(n-1)!}-\sum_{n=1}^{\infty}\frac{(1)_{n,\lambda}}{n!}t^{n} \nonumber \\
&=\sum_{n=1}^{\infty}n(1)_{n-1,\lambda}\frac{t^{n}}{n!}-\sum_{n=1}^{\infty}\frac{(1)_{n-1,\lambda}(1-(n-1)\lambda)}{n!}t^{n}\nonumber \\
&=\sum_{n=2}^{\infty}\big\{(n-1)(1)_{n-1,\lambda}+\lambda(n-1)(1)_{n-1,\lambda}\big\}\frac{t^{n}}{n!}\nonumber \\
&=\sum_{n=2}^{\infty}(1+\lambda)(n-1)(1)_{n-1,\lambda}\frac{t^{n}}{n!}. \nonumber
\end{align}
Therefore, by \eqref{70} and \eqref{71}, we obtain the following theorem. 
\begin{theorem}
For $n\ge 2$, we have 
\begin{displaymath}
\sum_{k=2}^{n}(1)_{k-1,\frac{1}{\lambda}}\lambda^{k-2}S_{2,\lambda}(n,k)=(n-1)(1)_{n-1,\lambda}. 
\end{displaymath}
\end{theorem}
We remark that, by letting $ \lambda \rightarrow 0$, we get the following corollary.
\begin{corollary}
For $ n \ge 2$, we have
\begin{align*}
\sum_{k=2}^{n}(-1)^{k-2}(k-2)!S_{2}(n,k)=n-1.
\end{align*}
\end{corollary}

\section{Conclusion}
In this paper, by exploiting generating functions we investigated some properties, recurrence relations and identities involving degenerate Stirling numbers associated with degenerate hyperharmonic numbers and also with degenerate Bernoulli, degenerate Euler, degenerate Bell and degenerate Fubini polynomials. In particular, the degenerate hyperharmonic numbers were introduced as a degenerate version of the hyperharmonic numbers. \par
It is one of our future projects to continue to study various degenerate versions of special polynomials and numbers by using such diverse methods as it is mentioned in Introduction.

\end{document}